\documentclass{article}

\usepackage{amsmath}
\usepackage{amsthm}
\usepackage{amssymb}
\usepackage{enumerate}
\usepackage[all]{xy}

\DeclareSymbolFont{AMSb}{U}{msb}{m}{n}
\DeclareMathSymbol{\N}{\mathbin}{AMSb}{"4E}
\DeclareMathSymbol{\Z}{\mathbin}{AMSb}{"5A}
\DeclareMathSymbol{\R}{\mathbin}{AMSb}{"52}
\DeclareMathSymbol{\Q}{\mathbin}{AMSb}{"51}
\DeclareMathSymbol{\I}{\mathbin}{AMSb}{"49}
\DeclareMathSymbol{\C}{\mathbin}{AMSb}{"43}

\newtheorem{theorem}{Theorem}[section]
\newtheorem{lemma}[theorem]{Lemma}
\newtheorem{prop}[theorem]{Proposition}
\newtheorem{corollary}[theorem]{Corollary}

\theoremstyle{definition}
\newtheorem{definition}[theorem]{Definition}
\newtheorem{remark}[theorem]{Remark}

\begin{document}

\title{A Composition Formula for Asymptotic Morphisms}         % Enter your title between curly braces
\author{J. Matthew Mahoney}        % Enter your name between curly braces
\date{}          % Enter your date or \today between curly braces
\maketitle
 
\begin{abstract}
 \noindent For graded $C^*$-algebras $A$ and $B$, we construct a semigroup ${\cal AP}(A,B)$ out of asymptotic pairs. This semigroup is similar to the semigroup $\Psi(A,B)$ of unbounded KK-modules defined by Baaj and Julg \cite{BaaJul} and there is a map $\Psi(A,B) \to {\cal AP}(A,B)$ when $B$ is stable. Furthermore, there is a natural semigroup homomorphism ${\cal AP}(A,B) \to E(A,B)$, where $E(A,B)$ is the E-theory group. We denote the image of this map $E'(A,B)$ and prove both that $E'(A,B)$ is a group and that the composition product of E-theory specializes to a composition product on these subgroups. Our main result is a formula for the composition product on $E'$ under certain operator-theoretic hypotheses about the asymptotic pairs being composed. This result is complementary to known results about the Kasparov product of unbounded KK-modules.
\end{abstract}

\section{Introduction}
The Kasparov product of unbounded KK-modules and the composition product of asymptotic morphisms in E-theory are two crucial structures in operator K-theory. Unfortunately for practitioners, at present, neither of these products is constructive in the sense that the data (unbounded KK-modules or asymptotic morphisms) determines a unique representative of the composition class. Nevertheless, there has been quite a bit of work on identifying the Kasparov product of unbounded KK-modules, culminating in the work of Kucerovsky \cite{Kuc} who proves general conditions under which a proposed KK-module is in fact the composition of two other modules. In this case, though, one does not have a formula for the composition, but rather a means of testing a given module to determine if it is the composition. 

For (graded) $C^*$-algebras $A$ and $B$, there is a canonical group homomorphism $KK(A,B) \to E(A,B)$. It is well-known that if $A$ is nuclear (or K-nuclear), then this map is an isomorphism and the two bivariant theories are the same. In the non-nuclear case, these groups are not necessarily isomorphic, but the image of a KK-element in E-theory induces the same map on K-theory that the KK-element does. In other words, at the level of K-theory there is no distinction between $KK(A,B)$ and its image in $E(A,B)$. 

In this paper we would like to compose unbounded KK-modules at the level of E-theory. In order to do this we must broaden the class of unbounded KK-modules to what we call {\it asymptotic pairs}. Our construction is similar to the construction in \cite{BaaJul} where Baaj and Julg build a semigroup $\Psi(A,B)$ out of unbounded KK-modules and proved that there is a surjection of semigroups $\Psi(A,B) \to KK(A,B)$. In our case, we prove that there is a semigroup ${\cal AP}(A,B)$ whose elements are asymptotic pairs and that there is a semigroup homomorphism ${\cal AP}(A,B) \to E(A,B)$ which is onto the image of $KK(A,B)$. The class of asymptotic pairs includes all unbounded KK-modules, but the relationship between the components of an unbounded KK-module is stronger than is needed to build an asymptotic morphism. Asymptotic pairs are built out of (almost) the same data as unbounded KK-modules, but the relationship between the components is weaker. 

The virtue of the weaker notion of asymptotic pair is that under certain operator-theoretic conditions a formula for the E-theory composition of asymptotic pairs is available. The proof of the composition formula crucially uses the weaker relationship. In particular, the representative of the composition given by the formula might not be an unbounded KK-module even if the two elements being composed are.

In the process of setting up the framework for the composition formula, we will see that the image, $E'(A,B)$, of the semigroup ${\cal AP}(A,B)$ in $E(A,B)$ is a group and is functorial in $A$ and $B$. By definition it is stable and homotopy invariant. Thus, we get (almost for free) that the composition product of E-theory specializes to a composition product on $E'$. In this way, we have a subcategory, $E'$, of $E$ which contains the image of the canonical functor $KK \to E$. The composition formula presented here is a special case of the general composition product in the category $E'$.

The paper is organized as follows. In Section \ref{SEC: AsyMor} we will discuss E-theory in the ``spectral picture'' of Guentner and Higson \cite{HigGue}. This is the picture best adapted to study the image of KK-theory. In Section \ref{SEC: AsyPair} we define asymptotic pairs and establish their relationship to unbounded KK-modules and E-theory. In Section \ref{SEC: Comp} we prove the composition formula for asymptotic pairs under appropriate operator-theoretic conditions. In Section \ref{SEC: Applications} we discuss two applications of the composition formula. There are two appendices. The first contains preliminary material on regular operators on Hilbert modules and unbounded multipliers of $C^*$-algebras. This material is well-known, but we use the opportunity to establish notation and pool the relevant facts. The second appedix contains a detailed computation that is required in the proof of the compostion formula, but which is long enough to detract from the flow of ideas.

\section*{Acknowledgements}
Much of the work in this paper appears in the author's PhD thesis at Dartmouth College. It is a pleasure to thank my advisor Jody Trout for his help throughout. I would also like to thank Dana Williams, Erik van Erp, and Dorin Dumitrascu for valuable conversations and reading preliminary drafts of this paper. Some of the results presented here rely on unpublished, preliminary work of Nigel Higson and Jody Trout. I would like to thank them for communicating their work to me.

\section{Asymptotic Morphisms and E-theory}\label{SEC: AsyMor}
Asymptotic morphisms are the basic building blocks of E-theory \cite{ConHig}. We recall the important definitions and constructions here.

\begin{definition} A family of maps $\{\phi_t\}_{t\in [1,\infty)}:A \rightarrow B$ is called an {\bf asymptotic morphism} if the following hold:

\begin{description}
 \item[AM1] For all $a \in A$, $\phi_t(a)$ is continuous in $t$.
 \item[AM2] For all $a,b,c \in A$ and $\lambda \in \C$,
  \begin{enumerate}
  \item[(i)] $\lim\limits_{t\rightarrow \infty}( \phi_t(a+\lambda b) -\phi_t(a)-\lambda \phi_t(b) ) = 0$ \\
  \item[(ii)] $\lim\limits_{t\rightarrow \infty}( \phi_t(ab)  - \phi_t(a)\phi_t(b) ) = 0$ \\
  \item[(iii)] $\lim\limits_{t\rightarrow \infty}( \phi_t(a^*)  -\phi_t(a)^* ) = 0$ \\
  \item[(iv)] $\lim\limits_{t \to \infty}( \phi_t(\gamma_A(a)) - \gamma_B(\phi_t(a)) = 0$.
  \end{enumerate}
\end{description}

If condition {\bf AM1} holds for a family of maps $\{\phi_t\}_{t\in [1,\infty)}:A \rightarrow B$, then we say that it is a {\bf continuous family}.
\end{definition}
An asymptotic morphism gives rise to an actual $*$-homomorphism, but into a new $C^*$-algebra.

\begin{definition}
 Given a $C^*$-algebra $A$, let $A_{\infty} = C_b([1,\infty),A) / C_0([1,\infty),A)$. We call $A_{\infty}$ the {\bf infinity algebra} of $A$. The image of $a \in A$ in $A_{\infty}$ is denoted $[a]$. If $A$ is graded by $\gamma$, then $C_b([1,\infty),A)$ has the grading given by applying $\gamma$ pointwise and $C_0([1,\infty),A)$ is a graded ideal. Thus, $A_{\infty}$ is graded as well.
\end{definition}

It is straightforward to show that an asymptotic morphism $\{\phi_t\}:A \to B$ defines a graded $*$-homomorphism $\phi:A \to B_{\infty}$ by $a \mapsto [\phi_t(a)]$. Two asymptotic morphisms are called asymptotically equivalent if they define the same $*$-homomorphism in to $B_{\infty}$. One of the important facts in E-theory is that, using the Bartle-Graves Selection theorem \cite{BarGra}, any $*$-homomorphism $\phi:A \to B_{\infty}$ can be continuously lifted to a map into $C_b([1,\infty),B)$, and hence yields an asymptotic morphism. The lift is not unique, but any two lifts are asymptotically equivalent. One of the common techniques for constructing an asymptotic morphism is to define a $*$-homomorphism into an infinity algebra. We will not distinguish between an asymptotic morphism and its associated $*$-homomorphism unless it is necessary to avoid confusion.

The operation that sends a $C^*$-algebra $A$ to its infinity algebra $A_{\infty}$ is an exact functor (Lemma 3.3 \cite{Gue2}). Recall that a functor from the category of (graded) $C^*$-algebras to itself is exact if it takes exact sequences to exact sequences. The exactness of taking infinity algebras has the following easy implication.

\begin{prop}
 Let $J$ be an ideal in $B$. If $\phi:A \to B_{\infty}$ is a $*$-homomorphism such that $\pi \circ \phi = 0$, where $\pi:B_{\infty} \to (B/J)_{\infty} \cong B_{\infty}/J_{\infty}$ is the canonical projection, then there is an asymptotic morphism $\phi':A \to J_{\infty}$ which is asymptotically equivalent to $\phi$.
\end{prop}

\begin{remark}
 If $\{\phi_t\}:A \to B$ is an asymptotic morphism and its associated homomorphism $\phi:A \to B_{\infty}$ is as in the above corollary, then we say that $\{\phi_t\}$ is asymptotically in $J$. When defining an asymptotic morphism into a $C^*$-algebra $J$, it suffices to define it by constructing an asymptotic morphism into any algebra containing $J$ as an ideal and showing that it is asymptotically in $J$. In effect, this enlarges the equivalence classes that comprise elements of $J_{\infty}$.
\end{remark}

The important equivalence relation on asymptotic morphisms is homotopy. Two asymptotic morphisms $\{\phi^0_t\}:A \to B$ and $\{\phi^1_t\}:A \to B$ are homotopic if there exists an asymptotic morphism $\{H_t\}:A \to B[0,1]$ such that the the composition of $H$ with the evaluation homomorphisms at 0 and 1 yield $\phi^0$ and $\phi^1$ repectively. The collection of homotopy classes of asymptotic morphisms between $A$ and $B$ is denoted $[\![A,B]\!]$ and the homotopy class of an asymptotic morphism $\{\phi_t\}$ is denoted $[\![\phi_t]\!]$.

In \cite{ConHig}, Connes and Higson proved that there is an associative pairing on the set of homotopy classes of asymptotic morphisms which extends the usual composition of (homotopy classes of) $*$-homomorphisms. Namely, there is a map $$\circ:[\![A,B]\!] \times [\![B,C]\!] \to [\![A,C]\!]$$ which is simply composition if the classes represented are represented by $*$-homomorphisms and is associative under iterated pairings. This is nontrivial because the naive composition of asymptotic morphisms as families of maps is not, in general, an asymptotic morphism (e.g., $\psi_t\circ\phi_t$ need not satisfy any of the limit conditions). The construction of the product boils down to finding an increasing rescaling function $r:[1,\infty) \to [1,\infty)$ such that $\psi_{r(t)} \circ \phi_t$ is an asymptotic morphism. The remainder of the proof requires showing that this is well-defined at the level of homotopy. 

The E-theory of Connes and Higson is constructed out of homotopy classes of asymptotic morphisms. In the category of graded $C^*$-algebras Guentner and Higson \cite{HigGue} use some extra structure to give a spectral picture of E-theory. In this picture, the E-theory group is defined as follows $$E(A,B)=[\![S\hat{\otimes}A,B\hat{\otimes}{\cal K}]\!],$$ where $S$ denotes the $C^*$-algebra $C_0(\R)$ graded into even and odd functions (c.f. Appendix \ref{SEC: RegOps}). Thus, the E-theory group $E(A,B)$ is the set of homotopy classes of asymptotic morphisms from $S\hat{\otimes}A$ to $B\hat{\otimes}{\cal K}$. The group structure is given by direct sum using the fact that $M_2({\cal K}) \cong {\cal K}$. See \cite{HigGue} for details of the construction. 

In the category of graded $C^*$-algebras there is a ``comultiplication'' $\Delta$ given by
$$\xymatrix@R=.05in{
 S \ar[r]^-{\Delta} & S\hat{\otimes}S \\
 f \ar@{|->}[r] & f(X\hat{\otimes}1 + 1\hat{\otimes}X).}$$
This map is simply the functional calculus of the unbounded self-adjoint multiplier $X\hat{\otimes}1 + 1\hat{\otimes}X$ on $S\hat{\otimes}S,$ where $X$ is the self-adjoint multiplier $f(x) \mapsto xf(x)$ on $S$. We omit the details here, but it is easy to show that $X\hat{\otimes}1 + 1\hat{\otimes}X$ is in fact essentially self-adjoint on the dense, graded $*$-subalgebra $C_c(\R) \hat{\odot} C_c(\R)$. It is clear that $X \hat{\otimes}1$ and $1 \hat{\otimes}X$ commute and Corollary \ref{COR: fclemma2} (below) allows us to calculate this map explicitly.

There is a composition product on E-theory given by the composition of asymptotic morphisms. Specifically, if $[\![\phi_t]\!] \in E(A,B)$ and $[\![\psi_t]\!] \in E(B,C)$ then the product in E-theory is given by the following composition of homotopy classes asymptotic morphisms

$$\xymatrix{S\hat{\otimes}A \ar[r]^-{[\![\Delta\hat{\otimes}1]\!]} &S\hat{\otimes}S\hat{\otimes}A \ar[r]^-{[\![1\hat{\otimes}\phi_t]\!]} &S\hat{\otimes}B \ar[r]^-{[\![\psi_t]\!]} &C.}$$ The arrows in this diagram are arrows in the homotopy category of asymptotic morphisms. The presence of the comultiplication in the composition is a feature of the fact that we are working with graded $C^*$-algebras. 

\begin{remark}\label{REM: gradKtheory}
 The reason this is called the ``spectral picture'' of E-theory is because Trout \cite{Tro} shows that in the case $A = \C$ the group $E(\C,B) \cong \hat{K}(B)$, where $\hat{K}(B)$ is the graded K-theory group defined by Kasparov. The group $\hat{K}(B)$ is built out of unbounded operators on Hilbert $B$-modules and the map $\hat{K}(B) \to E(\C,B)$ is given by sending an unbounded operator to the class of its functional calculus homomorphism. 
\end{remark}

\section{Asymptotic Pairs}\label{SEC: AsyPair}
We now define asymptotic pairs and prove their relationship to unbounded KK-modules. Throughout this section $A$ and $B$ will be graded $C^*$-algebras. To motivate the definition of asymptotic pairs we recall the definition of unbounded KK-modules which originally appeared in \cite{BaaJul}.

\begin{definition}
 An {\bf unbounded KK-module} for $(A,B)$ is a triple $(X_B,\phi,D)$ where $X_B$ is a countably generated, graded right Hilbert $B$-module, $\phi:A \to {\cal L}(X_B)$ is a nondegenerate graded $*$-homomorphism, and $D$ is a regular, odd unbounded self-adjoint operator on $X_B$ satisfying the following conditions:
\begin{description}
 \item[KK1]\label{KK1} For any $f \in S$ and $a \in A$, $f(D)\phi(a) \in {\cal K}(X_B)$.
 \item[KK2]\label{KK2} The set of $a \in A$, such that $[D,a]$ extends to an element of ${\cal L}(X_B)$ is dense in $A$.
\end{description}
The set of unbounded KK-modules for $(A,B)$ is denoted by $\Psi(A,B)$.
\end{definition}

\begin{remark}
 The condition that $\phi$ be nondegenerate is not in the original definition of an unbounded KK-module. However, one can assume that this is true up to KK-equivalence (c.f. Proposition 18.3.6 in \cite{Bla}). Note that what is called an essential homomorphism in \cite{Bla}, we call a nondegenerate homomorphism.
\end{remark}

In \cite{BaaJul}, Baaj and Julg proved that when $B$ is a stable $C^*$-algebra $\Psi(A,B)$ can be made into a semigroup and that there is a natural surjective semigroup homomorphism $\Psi(A,B) \to KK(A,B)$. The following definition is closely modeled on the the notion of an unbounded KK-module, but is better adapted to working with asymptotic morphisms.

\begin{definition}
 Let $\phi:A \to M(B)$ be a strict $*$-homomorphism and $D$ be an odd, unbounded self-adjoint multiplier of $B$. We say that $(\phi,D)$ is an {\bf asymptotic pair} if for every $f \in S$ and $a \in A$ we have
\begin{description}
 \item[AP1] \label{AP1} $f(D)\phi(a) \in B$ and
 \item[AP2] \label{AP2} $[f(t^{-1}D),\phi(a)] \to 0$ as $t\to \infty$.
\end{description}
\end{definition}
The following proposition constructs an asymptotic morphism from an asymptotic pair. In the case that $B$ is stable, this construction provides a map from ${\cal AP}(A,B)$ to $E(A,B)$.

\begin{prop}
 Let $(\phi,D) \in {\cal AP}(A,B)$. The assignment $$f \hat{\otimes} a \mapsto f(t^{-1}D)\phi(a)$$ defines up to equivalence an asymptotic morphism from $S \hat{\otimes} A$ to $B$.
\end{prop}
\proof
$f \mapsto f(t^{-1}D)$ is a continuous family of $*$-homomorphisms from $S \to M(B)$. By {\bf AP2} we know that the continuous family of $*$-homomorphisms $f(t^{-1}D)$ asymptotically commutes with the $*$-homomorphism $\phi$. Thus, by the universal property of the maximal graded tensor product, the product $f(t^{-1}D)\phi(a)$ defines a map from $S\hat{\otimes}A$ to $B_{\infty}$, hence an asymptotic morphism.
\endproof

We will now define a semigroup structure on ${\cal AP}(A,B\hat{\otimes}{\cal K})$. For notational convenience, we will assume that $B$ is a stable $C^*$-algebra and we have fixed an even $*$-isomorphism $\beta:B \to M_2(B)$. There is only one such graded isomorphism at the level of homotopy.

We define a semigroup structure on the set of asymptotic pairs as follows. Let $(\phi,D), (\tilde{\phi},\tilde{D})  \in {\cal AP}(A,B)$, then their sum is by definition $$(\phi,D)+(\tilde{\phi},\tilde{D}) = ({\rm diag}(\phi,\tilde{\phi}), {\rm diag}(D,\tilde{D})),$$ where we identify elements in $M_2(B)$ with $B$ through the isomorphism $\beta$. The semigroup structure depends on the choice of isomorphism $\beta$, but for any choice of $\beta$, the canonical map ${\cal AP}(A,B) \to E(A,B)$ is a semigroup homomorphism and the image of ${\cal AP}(A,B)$ is independent of $\beta$. We denote the image of this map by $E'(A,B)$. The following proposition states that $E'(A,B)$ is in fact a group.

\begin{prop}
\label{PROP: group}
 The additive inverse of $(\phi,D)$ is $(\phi^{{\rm opp}},-D)$, where $\phi^{\rm opp}$ is the opposite homomorphism of $\phi$.
\end{prop}
\proof
From Lemma 2.1 in \cite{HigGue}, the only thing we need to prove is that the class of asymptotic pairs is ``closed under opp''. Let $\Phi_t$ be the asymptotic morphism associated to $(\phi,D)$. The induced grading on $S\hat{\otimes}A$ is given by $f(x)\hat{\otimes}a \mapsto f(-x)\hat{\otimes}\gamma(a)$. It is clear, then, that $\Phi_t^{{\rm opp}}$ is the asymptotic morphism associated to $(\phi^{{\rm opp}},-D)$.
\endproof

From Propoition \ref{PROP: group}, we have that $E'(A,B)$ is a group. We will now show that there exists a map $\Psi(A,B) \to {\cal AP}(A,B)$ when $B$ is a stable $C^*$-algebra. 

\begin{prop}
 Let $A$ and $B$ be graded $C^*$-algebras and let $(X_B,\phi,D) \in \Psi(A,B)$.  Then $(\phi,D) \in {\cal AP}(A,{\cal K}(X_B))$. If $B$ is a stable $C^*$-algebra, then $(\phi,D) \in {\cal AP}(A,B)$.
\end{prop}
\proof

The proof that {\bf AP1} holds is immediate because $\phi$ is nondegenerate by assumption. To demonstrate {\bf AP2}, observe that by definition there exists a dense, graded $*$-subalgebra ${\cal A} \subseteq A$ such that the commutator $[D,a]$ is a bounded operator. Furthermore, we know that the $*$-algebra ${\cal S}$ generated by the resolvent functions $r_{\pm}(x)=(x\pm i)^{-1}$ is dense in $S$. 

Let $a \in {\cal A}_0$. The resolvent identity $$[r_{\pm}(t^{-1}D),a] = -r_{\pm}(t^{-1}D)[t^{-1}D,a]r_{\pm}(t^{-1}D)$$ yields the estimate $$\|[r_{\pm}(t^{-1}D),a]\| \leq t^{-1}\|[D,a]\|.$$ For the case where $a \in {\cal A}_1$, first split the resolvent into even and odd parts and collect terms so that the (graded) commutator of $D$ and $a$ appears in the middle of each term.

Thus, $[r_{\pm}(t^{-1}D),a] \to 0$ as $t \to \infty$ for any $a \in {\cal A}$. Since ${\cal S}$ is dense in $S$ and ${\cal A}$ is dense in $A$, an approximation argument finishes the proof of the first assertion.

To prove the second assertion, note that if $B$ is stable, then there is a canonical isomorphism $B \cong {\cal K}(\hat{{\cal H}}_B)$, where $\hat{{\cal H}}_B$ is the standard graded Hilbert $B$-module. By the Kasparov Stabilization theorem \cite{Kas}, there is a degree zero unitary isomorphism $X_B \oplus \hat{{\cal H}}_B \cong \hat{{\cal H}}_B$. Thus, we can include ${\cal K}(X_B)$ as a corner in ${\cal K}(\hat{{\cal H}}_B) \cong B$ and any two such identifications are homotopic by a graded version of  Higson's generalization of Kuiper's theorem \cite{Hig}. The second assertion follows.
\endproof

Thus, there is a commutative diagram of set maps $$\xymatrix{\Psi(A,B) \ar[r] \ar[d] & {\cal AP}(A,B) \ar[d]\\ KK(A,B) \ar[r] & E(A,B)}$$ where the vertical maps are semigroup homomorphisms and the bottom map is a homomorphism of groups. Since the canonical map $\Psi(A,B) \to KK(A,B)$ is surjective, the map ${\cal AP}(A,B) \to E(A,B)$ is onto the image of $KK(A,B)$. The semigroup ${\cal AP}(A,B)$ is {\em a priori} larger than $\Psi(A,B)$. We will make explicit use of this when composing asymptotic pairs.

The group homomorphism $KK(A,B) \to E'(A,B)$ shows that $E'(A,B)$ is a bivariant functor from the category of (separable) graded $C^*$-algebras with $*$-homomorphisms to the category of abelian groups. Thus $E'(\cdot,\cdot)$ is a stable, additive, homotopy invariant functor. From (a graded version of) Higson's characterization of E-theory as a category of fractions, we know that we have a pairing $$E'(A,B) \times E(B,C) \to E'(A,C),$$ which extends the usual composition of $*$-homomorphisms (c.f. Theorem 25.6.1 in \cite{Bla}). Since $E'(B,C)$ is a subgroup of $E(B,C)$, this pairing restricts to a pairing $$E'(A,B) \times E'(B,C) \to E'(A,C).$$  In this way we have constructed a full additive subcategory $E'$ of $E$. In particular, the composition of asymptotic pairs is an asymptotic pair.

\begin{remark}
 The above result holds because we declare two asymptotic pairs to be equivalent if their asymptotic morphisms are equivalent in $E(A,B)$. It is tempting to declare two asymptotic pairs to be equivalent if there is an ``asymptotic pair homotopy'' between them (where an asymptotic pair homotopy is an asymptotic pair that represents an E-theory equivalence between the pairs). This definition suffers because the homotopy in Lemma 2.1 of \cite{HigGue} is not obviously from an asymptotic pair and we would (potentially) lose the group structure on the equivalence classes of asymptotic pairs.
\end{remark}

\begin{remark}
 The definition of an asymptotic pair suggests that the subgroup $E'(A,B)$ of $E(A,B)$ is simply the image of $KK(A,B)$ under the canonical homomorphism. This question is still open, but it is not {\em a priori} clear either that $\Psi(A,B) = {\cal AP}(A,B)$ or that the equivalence relation placed on KK-modules by the map $$\Psi(A,B) \to {\cal AP}(A,B) \to E(A,B)$$ is the same relation placed on KK-modules by the map $$\Psi(A,B) \to KK(A,B).$$ Asymptotic pairs are a generalization of unbounded KK-modules and the natural (and weaker) equivalence relation to place on them is ``asymptotic homotopy equivalence''. 
\end{remark}
                                                                                                                                                                                                
\section{Composition}\label{SEC: Comp}
From the previous section, we know that the composition product on E-theory specializes to a composition product on $E'$. In this section we provide a formula for the product under certain operator-theoretic assumptions about the pairs to be composed. To motivate the following results, let $(\phi,D) \in {\cal AP}(A,B)$ and $(\psi,D') \in {\cal AP}(B,C)$ and suppose that $\psi(D)$ and $D'$ commute as odd, unbounded self-adjoint multipliers of $C$, i.e. that $[\psi(D),D']=\psi(D)D'+D'\psi(D)=0$. (The commutator of two unbounded multipliers is a subtle thing \cite{Kuc} and we will only treat this commutator formally for now.) For notational convenience, we will identify $D$ with its push-forward $\psi(D)$ and we will identify elements of $A$ and $B$ with their respective images under $\phi$, $\psi$, and $\psi \circ \phi$. Let $\phi_t$ and $\psi_t$ be the asymptotic morphisms associated to $(\phi,D)$ and $(\psi,D')$ respectively. We begin with elementary tensors of the form $e^{-x^2} \hat{\otimes}a$. Naively, we expect
$$\xymatrix@R = .1in{
S\hat{\otimes}A \ar[r]^-{\Delta\hat{\otimes}1} &S\hat{\otimes}S\hat{\otimes}A \ar[r]^-{1\hat{\otimes}\phi_t} &S\hat{\otimes}B \ar[r]^-{\psi_t} &C. \\ e^{-x^2}\hat{\otimes}a \ar@{|->}[r] & e^{-x^2}\hat{\otimes}e^{-x^2}\hat{\otimes}a \ar@{|->}[r] & e^{-x^2}\hat{\otimes}e^{-t^{-2}D^2}a \ar@{|->}[r] & e^{-t^{-2}D'^2}e^{-t^{-2}D^2}a
}$$
to be the compositon of $\psi_t$ and $\psi_t$ on tensors of the form $e^{-x^2} \hat{\otimes}a$.

Note that $e^{-t^{-2}D^2}a \in B$ because $(\phi,D) \in {\cal AP}(A,B)$. Now, formally at least, we have that $$e^{-t^{-2}(D'+D)^2} = e^{-t^{-2}(D'^2+D'D+DD'+D^2)} = e^{-t^{-2}D'^2}e^{-t^{-2}D^2}$$ because the cross terms in the exponent cancel and $D^2$ and $D'^2$ are even multipliers which commute with each other. A similiar computation with the function $h(x)=xe^{-x^2}$ shows that $h \hat{\otimes} a \mapsto h(t^{-1}(D+D'))a$ is (formally) equal to the naive composition. Since the algebra generated by these elementary tensors is dense, this suggests that the composition of $(\phi,D)$ and $(\psi,D')$ ought to be $(\psi \circ \phi,D+D')$.

The above computation is far from complete. For one, there are many arithmetic operations with unbounded multipliers. Indeed, the sum $D+D'$ may not even make sense because its natural domain ${\rm dom}(D) \cap {\rm dom}(D')$ may not be dense. Second, it need not be the case that the last map is well-defined as an asymptotic morphism. The remainder of the paper is focused on providing appropriate technical hypotheses for this computation to go through. In the process, we will only require that the commutator of $D$ and $D'$ be bounded (not zero), which is a substantial improvement over the above. The intuition behind why we can accept this much weaker hypothesis is that if $[D,D']$ is bounded, then $[t^{-1}D,t^{-1}D'] \to 0$ as $t\to \infty$, i.e. $D$ and $D'$ asymptotically commute. 

\begin{remark}
There was a choice made to require the hypothesis of strictness in the definition of an asymptotic pair. If that hypothesis is dropped, then the statement of the composition formula would require that it be added in for push-forwards to make sense. On the other hand, using the theory of connections, one can prove that an unbounded KK-module can be replaced up to KK-equivalence by one with a nondegenerate homomorphism. Thus, the more restrictive definition of asymptotic pair retains a map $KK(A,B) \to E'(A,B)$ which factors the canonical map from KK-theory to E-theory. 
\end{remark}

\subsection{Unbounded multipliers and Commutators}
In this subsection we collect the necessary technical hypotheses and estimates concerning the commutator of unbounded multipliers.

Let $N>0$ and let ${\rm i}_N(x)=x(1+N^{-2}x^2)^{-1} \in C_b(\R)$. It is clear that ${\rm i}_N(x)$ approaches the identity function $x \mapsto x$ uniformly on compact subsets of $\R$ as $N\rightarrow \infty$. If $D$ is an odd, essentially self-adjoint multiplier of a graded $C^*$-algebra $C$, then we would expect $D_N = {\rm i}_N(D)$ to approximate $D$ somehow. Indeed, we have the following easy lemma.

\begin{lemma} \label{LEM: approxlem} Let $D$ be an odd, essentially self-adjoint multiplier of a $C^*$-algebra $C$ and $f \in C_0(\R)$. Then $f(D_N) \rightarrow f(D)$ in norm as $N \rightarrow \infty$.
\end{lemma}
\proof
These estimates hold for functions and therefore hold for the functional calculus.
\endproof

\begin{remark} If $f \notin C_c(\R)$, then $f \circ {\rm i}_N \notin C_0(\R)$. But $f \circ {\rm i}_N$ becomes arbitrarily small off of a compact set with increasing $N$. 
\end{remark}

In the composition formula and its proof, we will need to add unbounded multipliers to each other and make sense of their commutators. The following definition gathers the necessary technical assumptions together. 

\begin{definition} Let $C$ be a graded $C^*$-algebra and let $D$ and $D'$ be odd, essentially self-adjoint  multipliers of $C$. We say that $D$ and $D'$ have {\bf bounded commutator} if the following hold:
\begin{description}
\item[\bf{BC1}] The operators $D$, $D'$, and $D+D'$ have a common core $\mathcal{C} \subseteq C$.
\item[\bf{BC2}] The core $\mathcal{C}$ is mapped to itself by $D_N$ and $D'_N$ and isomorphically onto itself by the resolvents $(D \pm i)^{-1}$ and $(D' \pm i)^{-1}$.
\item[\bf{BC3}] The commutator $[D,D']$, defined on $\mathcal{C}$, extends to a bounded operator on $C$.
\end{description}
\end{definition}

\begin{remark} We include {\bf BC1} as one of the bounded commutator conditions because the reason we need to consider the commutator $[D,D']$ is to relate the functional calculus of $D+D'$ to those of $D$ and $D'$ separately. The technical conditions {\bf BC} are sufficient to prove that such a relationship exists.
\end{remark}

\begin{lemma} \label{LEM: commbound} Let $C$ be a graded $C^*$-algebra and let $D$ and $D'$ be odd, essentially self-adjoint  multipliers of $C$ which have bounded commutator. Then $$\| [D_N,D'_N] \| \leq \| [D,D'] \|,$$ for all $N > 0$.
\end{lemma}
 \proof
 The operator $D_N$ admits an integral decomposition (c.f. Proposition 2.2 in \cite{BaaJul})
 $$D_N =  \frac{2}{\pi} \int_0^{\infty} (N^{-2} D^2 + 1 + s^2)^{-1} D \ ds.$$
 \noindent This can be seen from the functional calculus because the relevant Riemann sums form a sequence in $C_0(\R)$ which converges uniformly on compacts subets of $\R$. Thus, the images of the Riemann sums in $M(C)$ form a strictly convergent sequence. For any $c \in \mathcal{C}$ we have that $D' c \in {\cal C}$ and $D_N c \in \mathcal{C}$ from {\bf BC2}. Thus,
 $$[D_N,D'] c = \frac{2}{\pi} \int_0^{\infty} [(N^{-2} D^2 + 1 + s^2)^{-1} D ,D'] c \  ds.$$ The commutator passes through the integral because $D'$ is closeable and the Riemann sums used to define the integral form a Cauchy sequence in $\mathcal{C}$ which converges to an element in $\mathcal{C}$. 
 
 Letting $p(x) = (N^{-1}x+i\sqrt{1+s^2})^{-1}$ and $q(x) = x p(x)\bar{p}(x)$ we have the identity of functions $$q(x) = \frac{N}{2} \big( p(x) + \bar{p}(x) \big)$$ and from the functional calculus we get the corresponding identity of operators $$q(D) = \frac{N}{2} \big( p(D) + \bar{p}(D) \big).$$ Note that $\bar{p}(D) = p(D)^*$. Because $q$ is an odd, bounded real-valued function, $q(D)$ is an odd, self-adjoint multiplier of $C$. Thus the (graded) commutator in the integrand is equal to $q(D) D' + D' q(D)$. Computing 
 \begin{align*}
 q(D) D' + D' q(D)  &= \frac{N}{2} \big(p(D) D' + D' \bar{p}(D) +  \bar{p}(D) D' + D' p(D)\big) \\ &= \frac{1}{2} \big( p(D) [D,D'] \bar{p}(D) + \bar{p}(D) [D,D'] p(D) \big).
 \end{align*}
 Note that there are no issues with the above factorizations because we are plugging in $c \in \mathcal{C}$, which will be mapped into the relevant domains by virtue of $\bf{BC2}$. From the functional calculus $\|p(D)\| \leq (1+s^2)^{-\frac{1}{2}}$. Thus we get the estimate
 $$ \| [q(D),D']  \| \leq (1+s^2)^{-1} \| [D,D'] \|.$$ It follows that $\| [D_N, D'] \| \leq \| [D,D'] \|$. To finish, perform the same calculation with $D'$ to get $$\|[D_N,D'_N]\| \leq \|[D_N,D']\|.$$ Note that this time there are no issues with domains because $D_N$ is bounded.
\endproof

\begin{corollary} If we let $D_t=t^{-1}D$ and $D_t'=t^{-1}D'$, then $\| [D_{t,N},D_{t,N}'] \| \leq \| [D_t,D_t'] \| $, where $D_{t,N}=(D_t)_N$ and $D_{t,N}'=(D_t')_N$. In particular, $\| [D_{t,N},D_{t,N}'] \| \rightarrow 0$ uniformly in $N$ as $t \rightarrow \infty$.
\end{corollary}

With the above lemma we can establish a crucial estimate that will relate the functional calculus of $D+D'$ to the calculi for $D$ and $D'$.

\begin{prop} \label{PROP: techlemma} For every $f \in C_0(\R)$, $$\lim_{N\rightarrow \infty} \limsup_{t \rightarrow \infty} \|f(D_{t,N}+D_{t,N}') - f(D_t+D_t') \| = 0.$$
\end{prop}
\proof
From the Stone-Weierstrass theorem, it suffices to prove this for the resolvent function $r(x) = (x+ i)^{-1}$. Consider the following tentative factorization:
\begin{align*}
r(D_{t,N}+D_{t,N}')\  -\  & r(D_t+D_t') \\&= r(D_{t,N}+D_{t,N}')(D_t - D_{t,N}) r(D_t+D_t') \\&\hspace{.5in} +  r(D_{t,N}+D_{t,N}')(D_t' - D_{t,N}') r(D_t+D_t')  \\ &= r(D_{t,N}+D_{t,N}')(1 - (1+N^{-2}D_t^2)^{-\frac{1}{2}})D_t r(D_t+D_t') \\&\hspace{.5in} +  r(D_{t,N}+D_{t,N}')(1 - (1+N^{-2}D_t'^2)^{-\frac{1}{2}})D_t' r(D_t+D_t').
\end{align*}
To rigorously establish this identity, we need to show that $D_tr(D_t+D_t')$ and $D_t'r(D_t+D_t')$ are bounded operators (the other factors on the right hand side are bounded). Since $D_t+D_t'$ is essentially self-adjoint on $\mathcal{C}$ by {\bf BC1}, we know that $D_tr(D_t+D_t')$ is densely defined. For $c \in (D_t+D_t' +i)\mathcal{C}$ we want to show that $$\|D_tr(D_t+D_t')c\|^2 \leq M \|c\|^2$$ for some constant $M>0$. From the $C^*$-identity we compute 
\begin{align*}
\|D_tr(D_t+D_t')c\|^2 &= \|c^*\bar{r}(D_t+D_t')D_t^2r(D_t+D_t')c\| \\ &\leq \|c^*\bar{r}(D_t+D_t')(D_t^2+D_t'^2)r(D_t+D_t')c\| \\ &\leq \|c^*\bar{r}(D_t+D_t')(D_t+D_t')^2r(D_t+D_t')c\| \\& \hspace{.5in}+ \|c^*\bar{r}(D_t+D_t')[D_t,D_t']r(D_t+D_t')c\| \\ &\leq \|c\|^2 + \|[D,D']\| \|c\|^2.
\end{align*}
Thus we may take $M = 1+\|[D,D']\|$. Similarly for the term with $D$ and $D'$ switched.

Now we want to show that the factor $$r(D_{t,N}+D_{t,N}')(1 - (1+N^{-2}D_t^2)^{-\frac{1}{2}}) \rightarrow 0$$ as $N \rightarrow \infty$, uniformly in $t$. From the $C^*$-identity and the functional calculus, it suffices to show that $$(1 + (D_{t,N}+D_{t,N}')^2)^{-1}(1 - (1+N^{-2}D_t^2)^{-\frac{1}{2}}) \rightarrow 0$$ as $N \rightarrow \infty$, uniformly in $t$. From Lemma \ref{LEM: commbound}, we know that $[D_{t,N},D_{t,N}'] \leq t^{-2}[D,D']$ for any $N$. Thus the difference $$\big(D_{t,N}+D_{t,N}'\big)^2 - \big(D_{t,N}^2+D_{t,N}'^2\big) = [D_{t,N},D_{t,N}'] \rightarrow 0$$ as $t \rightarrow \infty$. Since both $D_{t,N}^2$ and $D_{t,N}'^2$ are positive, $1+ D_{t,N}^2+D_{t,N}'^2$ is invertible. From the continuity of inversion, the difference $$(1 + (D_{t,N}+D_{t,N}')^2)^{-1} - (1+ D_{t,N}^2+D_{t,N}'^2)^{-1} \rightarrow 0$$ as $t \rightarrow \infty$. From positivity we have that $$(1+ D_{t,N}^2+D_{t,N}'^2)^{-1} \leq (1+ D_{t,N}^2)^{-1}.$$ Thus, it suffices to show that $$(1+ D_{t,N}^2)^{-1} (1 - (1+N^{-2}D_t^2)^{-\frac{1}{2}}) \rightarrow 0$$ as $N \rightarrow \infty$, uniformly in $t$. But this expression only involves the functional calculus for $D$ and a calculation with functions shows that the desired limit holds. The same argument works for the term with $D$ and $D'$ switched and the result follows.
\endproof

The above proposition has the following corollaries which will be used extensively in the calculations in the next section. Their proofs are involved and, therefore, deferred to the second appendix.

\begin{corollary} \label{COR: fclemma} Let $C$ be a graded $C^*$-algebra and let $D$ and $D'$ be odd, essentially self-adjoint multipliers of $C$ which have bounded commutator. Then $$\lim_{t \rightarrow \infty} \| e^{-t^{-2}(D+D')^2} -e^{-t^{-2}D^2}e^{-t^{-2}D'^2} \| = 0.$$ Furthermore, $$\lim_{t \rightarrow \infty} \big\| t^{-1}(D+D')e^{-t^{-2}(D+D')^2} - \big(t^{-1}De^{-t^{-2}D^2}\big)e^{-t^{-2}D'^2} - e^{-t^{-2}D^2}\big(t^{-1}D'e^{-t^{-2}D'^2}\big) \big\| = 0.$$
\end{corollary}

\begin{corollary} \label{COR: fclemma2}If $D$ and $D'$ are as above and $[D,D']=0$, then $$e^{-(D+D')^2} = e^{-D^2}e^{-D'^2}.$$ Furthermore, $$(D+D')e^{-(D+D')^2} = De^{-D^2}e^{-D'^2} + e^{-D^2}D'e^{-D'^2}.$$
\end{corollary}

\begin{definition}
Let $\phi_t:A \to C$ be an asymptotic morphism. We say that $\phi_t$ {\bf asymptotically factors} if there exist two asymptotic morphisms $\alpha_t:A \to B$ and $\beta_t:B \to C$ such that $\phi_t$ is asymptotically equivalent to the continuous family $\beta_t \circ \alpha_t$ (in which case, $\beta_t \circ \alpha_t$ is actually an asymptotic morphism).
\end{definition}

The next proposition summarizes the above work by saying that if $D$ and $D'$ have bounded commutator, then the functional calculus homomorphism of $t^{-1}(D+D')$ asymptotically factors.

\begin{prop}\label{PROP: asycomm}
 If $D$ and $D'$ are odd, unbounded multipliers of a $C^*$-algebra $C$ and have bounded commutator, then the continuous family of $*$-homomorphisms given by the functional calculus for $t^{-1}(D+D')$ asymptotically factors into the composition $$\xymatrix{S \ar[r]^-{\Delta} &S\hat{\otimes}S \ar[r]^-{\Phi_{D,D'}} & M(C)},$$ where the first map is the comultiplication $\Delta$ and the second map $\Phi_{D,D'}$ is the tensor product of the functional calculi for $t^{-1}D$ and $t^{-1}D'$.
\end{prop}
\proof
First we need to prove that $\Phi_{D,D'}:S\hat{\otimes}S \to M(C)$ given by $f \hat{\otimes}g \mapsto f(t^{-1}D)g(t^{-1}D')$ is in fact an asymptotic morphism. Note that it suffices to show that the functional calculus of $t^{-1}D$ and $t^{-1}D'$ asymptotically commute on  the generators $e^{-x^2}$ and $xe^{-x^2}$. Now, by Corollary \ref{COR: fclemma}, we know that $$e^{t^{-2}(D+D')^2} - e^{t^{-2}D^2}e^{t^{-2}D'^2} \to 0,$$ as $t \to \infty$. Since addition is commutative, we can asymptically commute $e^{t^{-2}D^2}$ and $e^{t^{-2}D'^2}$. For elementary tensors with $xe^{-x^2}$ in them, note that for large enough $N > 0 $, the difference $$xe^{-x^2} - {\rm i}_N(x)e^{-x^2}$$ can be made arbitrarily small. Since $D$ and $D'$ have bounded commutator, by Corollary \ref{COR: asycomm2} (below), we can commute any function of $D$ with $D_N'$ and vice versa. Thus, $\Phi_{D,D'}$ is an asymptotic morphism.

Since $\Delta$ is a $*$-homomorphism, the compositon of $\Delta$ and $\Phi_{D,D'}$ is straightforward (no rescaling is necessary). Corollary \ref{COR: fclemma} now gives us that $e^{-t^{-2}(D+D')^2}$ is asymptotic to $e^{-t^{-2}D^2}e^{-t^{-2}D'^2} = \Phi_{D,D'}(\Delta(e^{-x^2}))$. Similarly for the image of $xe^{-x^2}$. The result follows.
\endproof

\subsection{Composition Formula}
We now turn to a proof of the composition formula. The next proposition concerns the rigidity of asymptotic pairs. In general, by the exactness of the infinity algebra functor, one can define an asymptotic morphism into $B$ by constructing an asymptotic morphism into $M(B)$ which is asymptotically in $B$ (Remark \ref{REM: gradKtheory}). Such a family of maps would in principle never have to land in $B$. For asymptotic pairs, however, this is not the case.
                                                                                                                                                                                                
\begin{prop} \label{PROP: asyideal}
 Let $(\phi,D) \in {\cal AP}(A,B)$ and $J$ be an ideal in $B$. If the asymptotic morphism $f \hat{\otimes} a \mapsto f(t^{-1}D)\phi(a)$ is asymptotically in $J$, then $(\phi,D) \in {\cal AP}(A,J)$.
\end{prop}
\proof
Let $f \in C_c(\R)$ and suppose that the interval $[-N,N]$ contains the support of $f$. We wish to show that $f(D)\phi(a) \in J$. Let $\chi \in C_0(\R)$ be any function which is one on the interval $[-N,N]$. Then for all $t\in [1,\infty)$ we have the equality $$f(D)\phi(a) = \chi(t^{-1}D)f(D)\phi(a) = f(D)\chi(t^{-1}D)\phi(a).$$ The right hand side approaches $J$ as $t \to \infty$, but the left hand side is independent of $t$ and therefore must lie in $J$. The result follows because $C_c(\R)$ is dense in $C_0(\R)$.
\endproof    

The following lemma is a useful estimate. Its proof is a straightforward algebraic manipulation.

\begin{lemma}\label{LEM: gradedresolventidentities}
 Let $D$ be an odd, unbounded self-adjoint multiplier of $C$ and let $T \in M(C)$ be homogeneous. Suppose that the commutator $[D,T]$ is densely-defined and bounded. Then the following estimates hold $$\|[(D^2+1)^{-1},T]\| \leq \|[D,T]\|$$ and $$\|[D(D^2+1)^{-1},T]\| \leq \|[D,T]\|.$$
\end{lemma}
\proof
We will prove this for the commutator $[(D^2+1)^{-1},T]$ where $T$ is even. The rest are simple modifications.

By partial fractions we have the identity $(D^2+1)^{-1} = \frac{i}{2}\big((D+i)^{-1} - (D-i)^{-1}\big)$. Thus, the commutator splits into two terms. Since $(D^2+1)^{-1}$ is an even multiplier the commutator  is $$(D^2+1)^{-1}T - T(D^2+1)^{-1} = \frac{i}{2} \big( (D+i)^{-1} - (D-i)^{-1}\big)T - \frac{i}{2}T\big((D+i)^{-1} - (D-i)^{-1}\big).$$ Collecting the first and fourth term we have the identity 
\begin{align*}
 \frac{i}{2}\big((D+i)^{-1} T - T(D-i)^{-1}\big) &= \frac{i}{2}\big((D+i)^{-1}(T(D-i) - (D+i)T)(D-i)^{-1}\big) \\&= \frac{i}{2}(D+i)^{-1}[T,D](D-i)^{-1}.
\end{align*}
The second and third term can be similarly combined and the result follows.

The proof in the remaining cases requires the same trick, but combining terms in different ways to arrange that the graded commutator of $D$ and $T$ appears in the middle of each term.
\endproof

\begin{remark}
 The above result works by the same calculation for an unbounded multiplier $D'$ in place of $T$ provided that $D$ and $D'$ have bounded commutator.
\end{remark}

\begin{corollary}\label{COR: asycomm2}
 Let $D$ and $D'$ have bounded commutator. Then for any $f \in S$,  $[f(t^{-1}D'),D_N] \to 0$ as $t \to \infty$.
\end{corollary}
\proof
It suffices to prove this for the functions $f(x)=(x^2+1)^{-1}$ and $g(x)=x(x^2+1)^{-1}$. The proof in each of these cases is the same, and we only prove the result for $f$.

Since $D$ and $D'$ have bounded commutator, the expression $[D',D_N]$ is densely defined and bounded by the proof of Lemma \ref{PROP: techlemma}. From Lemma \ref{LEM: gradedresolventidentities} we have the estimate $$\|[f(t^{-1}D'),D_N]\| \leq \|[t^{-1}D',D_N]\| \to 0$$ as $t \to \infty$.
\endproof 

\begin{remark}
 Note that by the graded derivation property of commutators (Lemma \ref{LEM: gradder}), the above result can be extended to the commutator of $f(t^{-1}D')$ and any polynomial in $D_N$.
\end{remark}

The next lemma provides us with the final important piece to proving the composition formula.

\begin{lemma}
 Let $(\phi,D) \in {\cal AP}(A,B)$ and let $(\psi,D') \in {\cal AP}(B,C)$. Suppose that $D$ and $D'$ have bounded commutator. Then for any $f \in S$ and $a \in A$, we have $[f(t^{-1}D'),a] \to 0$ as $t \to \infty$.
\end{lemma}
\proof
Let $f \in S$ and $a \in A$ be homogeneous. We wish to show that $$[f(t^{-1}D'),a] = f(t^{-1}D')a - (-1)^{\partial f \partial a}af(t^{-1}D') \to 0$$ as $t \to \infty$. Consider the first term above $$f(t^{-1}D')a = f(t^{-1}D')(D_N^2+1)(D_N^2+1)^{-1}a.$$ From Lemma \ref{LEM: approxlem} we know that $(D_N^2+1)^{-1}a$ approches $(D^2+1)^{-1}a$ as $N \to \infty$. Since $(\phi,D)$ is an asymptotic pair, $(D^2+1)^{-1}a \in B$. Thus, we may choose an $N$ large enough so that $(D_N^2+1)^{-1}a$ is arbitrarily close to an element of $B$. Note that this element is homogeneous of the same degree as $a$. Since $f(t^{-1}D)$ asymptotically commutes with $(D^2+1)^{-1}a$, there is a $t_0$ large enough so that the commutator $[f(t^{-1}D),(D_N^2+1)^{-1}a]$ is smaller than any $\epsilon > 0$ provided $N$ is large enough and $t > t_0$.

To finish, observe that $f(t^{-1}D')$ asymptotically commutes with $(D_N^2+1)$ by Corollary \ref{COR: asycomm2}. Thus, $f(t^{-1}D')$ asymptotically commutes with both terms in the product and the result follows.
\endproof

\begin{remark}
 The need for the above lemma comes from the fact that $\phi(a)$ is in $M(B)$ and not in $B$. In the definition of asymptotic pairs, we make no assumption about the rate at which the commutator $[f(t^{-1}D'),b]$ goes to zero. In particular, the convergence is not uniform in $b$. We require such generality because in the case of an unbounded KK-module, the way we prove that $[f(t^{-1}D'),b] \to 0$ is by the estimate $\|[f(t^{-1}D'),b]\| \leq t^{-1}\|[D,b]\|$ (which holds for a dense set of functions $f$). But here we see that this estimate involves the norm of an unbounded derivation acting on $b$ and is not uniform in $b$! 
\end{remark}

As an easy corollary, we get the following.

\begin{corollary}\label{COR: hopscotch}
 Let $(\phi,D)$ and $(\psi,D')$ be as above. The $*$-homomorphism $\psi \circ \phi:A \to M(C)$ asymptotically commutes with $\Phi_{D,D'}$.
\end{corollary}                                                                                                                                                                           
                                                                                                                                                                                                
Thus we arrive at our main theorem, which puts the above results together.                                                                                                                      
                                                                                                                                                                                                
\begin{theorem}\label{THM: compform}                                                                                                                                                                              
 Let $(\phi,D) \in {\cal AP}(A,B)$ and $(\psi,D') \in {\cal AP}(B,C)$. Suppose that $\psi(D)$ and $D'$ have bounded commutator. Then $[\![\psi \circ \phi,\psi(D)+D']\!] = [\![\psi,D']\!] \circ [\![\phi,D]\!] \in E'(A,C)$. 
\end{theorem}
\proof
For notational convenience, we drop the symbol $\psi(D)$ and use $D$ for both the multiplier of $B$ and its push-forward to $C$. The functoriality of push-forwards (Proposition \ref{PROP: pushforward}) eliminates possible ambiguity. We also let $\rho = \psi \circ \phi$. From Proposition \ref{COR: hopscotch}, the asymptotic morphism $\Phi_{D,D'}':S\hat{\otimes}S \to M(C)$ asymptotically commutes with $\rho:A \to M(C)$. By Proposition \ref{PROP: asycomm} we know that the functional calculus of $t^{-1}(D+D')$ asymptotically factors as $\Phi_{D,D'} \circ \Delta$. Thus, $(\rho,D+D')$ is an asymptotic pair for $(A,M(C))$. 

The asymptotic morphism $S\hat{\otimes}S\hat{\otimes}A \to M(C)$ given by the tensor product of $\Phi_{D,D'}$ and $\rho$ is asymptotically in $C$, thus so is the composition $$\xymatrix{S \hat{\otimes}A \ar[r]^-{\Delta} &S\hat{\otimes}S\hat{\otimes}A \ar[r]^-{\Phi_{D,D'}} &C}.$$ Since the asymptotic morphism associated to $(\rho,D+D')$ asymptotically factors as the above composition, we know from Proposition \ref{PROP: asyideal} that $(\rho,D+D') \in {\cal AP}(A,C)$. 

We must now prove that $(\rho,D+D')$ is the composition of $(\phi,D)$ and $(\psi,D')$. But this follows because the naive composition is asymptotically equivalent to the asymptotic morphism associated to $(\rho,D+D')$. Since a continuous family of maps that is asymptotically equivalent to an asymptotic morphism is itself an asymptotic morphism, we are done.
\endproof

\begin{remark}
 In the above proof we actually prove something stronger than stated. The naive composition of asymptotic morphisms is not typically an asymptotic morphism, but, under the conditions of the theorem, the naive composition is an asymptotic morphism. Furthermore, the class of the naive composition can be identified as the class of $(\rho,D+D')$. Thus, in this case we have a composition of asymptotic classes, not just asymptotic homotopy classes. This is not too surprising, however. Given an asymptotic pair coming from a homomorphism or a KK-module, we often have to perform a nontrivial homotopy to construct its associated asymptotic pair (e.g. replacing a homomorphism with a strict homomorphism). 
\end{remark}

\section{Applications} \label{SEC: Applications}
In this section we will explore two applications of the above composition formula. The first is to the spectral proof of the Bott Periodicity theorem and the second is to the pertubation of K-homology classes of elliptic differential operators. Both of these results are well-known, but the content of this section is that both are naturally recognized as instances of the general composition formula. In essence, we show that the composition formula axiomatizes a standard computation in operator K-theory and exposes, in particular, the operator-theoretic underpinnings of such computations.

\subsection{Bott Periodicity}
We will not recapitulate a full proof of the Bott Periodicity theorem here and instead defer to Section 2 of \cite{HigGue}. What we aim to stress here is where in that line of reasoning we find an instance of the composition formula.

To establish notation let $V$ be a finite dimension Euclidean vector space and ${\rm Cliff_{\C}}(V))$ denote the complexified Clifford algebra of $V$. There is a grading on ${\rm Cliff_{\C}}(V))$ given by declaring even (odd) degree monomials in elements of $V$ to be even (odd). We let ${\mathcal C}_c^{\infty}(V) = C_c^{\infty}(V,{\rm Cliff_{\C}}(V))$ be the algebra of compactly supported ${\rm Cliff_{\C}}(V)$-valued functions on $V$ and ${\mathcal C}(V) = C_0(V,{\rm Cliff_{\C}}(V))$ be the $C^*$-algebra closure of ${\mathcal C}_c^{\infty}(V)$ of functions which vanish at infinity on $V$. The algebra ${\mathcal C}(V)$ is graded into functions which take values in the even/odd part of ${\rm Cliff_{\C}}(V)$. In addition, we let ${\mathcal H} = L^2(V,{\rm Cliff_{\C}}(V))$ be the space of square-integrable ${\rm Cliff_{\C}}(V)$-valued functions on $V$. The inner product on ${\mathcal H}$ is given by declaring fundamental monomials in an orthonormal basis of $V$ and $1$ to be orthogonal in ${\rm Cliff_{\C}}(V)$ and integrating over $V$. The resulting inner product is independent of choice of basis for $V$.

\begin{theorem}[Bott Periodicity]
For any finite dimensional vector space $V$, we have the following isomorphism $$K(\C) \cong K({\mathcal C}(V)).$$ 
\end{theorem}
\proof (sketch)

By Theorem 1.12 of \cite{HigGue}, it suffices to prove that there exists an element $b \in K({\mathcal C}(V)) \cong E'(\C,{\mathcal C}(V))$ and $\alpha \in E'({\mathcal C}(V),{\mathcal K}({\mathcal H}))$ such that the induced map $\alpha_*:K({\mathcal C}(V)) \to K(\C) \cong \Z$ maps $b$ to $1$. We define $b$ and $\alpha$ as asymptotic pairs. 

Let $\iota:\C \to M({\mathcal C}(V))$ denote the $*$-homomorphism $\lambda \mapsto \lambda I$. Clearly, $\iota$ is strict (it is unital). Let $C:{\mathcal C}_c(V) \to {\mathcal C}(V)$ denote the map $f(v) \mapsto v\cdot f(v)$. It is clear that $C$ is an odd, essentially self-adjoint multiplier of ${\mathcal C}(V)$ with a unique self-adjoint extension. Furthermore, an easy computation shows that $f(C) \in {\mathcal C}(V)$. Indeed, $f(C)$ is given by the map $g(v) \mapsto f(v) \cdot g(v)$ for every $f \in S$. Thus, the pair $(\iota,C)$ is an asymptotic pair. We let $b$ denote its class in $K({\mathcal C}(V))$.

Now we wish to construct an element $\alpha \in E'({\mathcal C}(V),{\mathcal K}({\mathcal H}))$. Let $\phi:{\mathcal C}(V) \to {\mathcal B}({\mathcal H}))$ be the (nondegenerate) representation of ${\mathcal C}(V)$ on ${\mathcal H}$ as multiplication operators. Let $D$ be the operator given as follows: $${\rm dom}(D) = C_c^{\infty}(V,{\rm Cliff}(V)),$$ and $$(Df)(v) = \sum_1^n \hat{e_i}(\frac{\partial f}{\partial x_i}(v)),$$ where $\hat{e}(g(v)) = (-1)^{\partial g}g(v)e_i$ denotes ($\pm$) right Clifford multiplication by $e$ for homogeneous $g$ and extended linearly. Clearly $D$ is odd and integration by parts shows that it is formally self-adjoint. This is the Dirac operator of $V$.

A standard argument in elliptic operator theory (see e.g. \cite{Gue1}) shows both that $D$ is essentially self-adjoint and that $f(D)\phi(g)$ is compact for every $f \in S$ and $g \in {\mathcal C}(V)$. Furthermore, the Leibniz rule gives that $\|[t^{-1}D,\phi(g)]\| \leq {\rm const}\  t^{-1}$ for smooth, compactly supported $g$. Thus, $[f(t^{-1}D),\phi(g)]$ extends to a bounded operator on ${\mathcal H}$ and goes to zero as $t$ goes to infinity. Hence, we have an asymptotic pair $(M,D) \in E'({\mathcal C}(V),\C)$. We let $\alpha$ denote its class in $E'({\mathcal C}(V),{\mathcal K}({\mathcal H}))$. Note that we are viewing $D$ now as a multiplier of ${\mathcal K}({\mathcal H})$, see Proposition \ref{PROP: hilbmodop}.

We now wish to compose $\alpha$ and $b$ using the composition formula. In order to do this, we must push $C$ forward to a multiplier of ${\mathcal K}({\mathcal H})$ and establish the bounded commutator conditions. 

The push forward of $C$ is is simply the operation of left Clifford multiplication by $v$, but this time thought of as an operator on ${\mathcal H}$. It is still odd and essentially self-adjoint on the set of compactly supported functions. (Note that $C$ is the same map on compactly supported functions, but is being realized as a multiplier on {\em two} different $C^*$-algebras, ${\mathcal C}(V)$ and ${\mathcal K}({\mathcal H})$.) 

In order to establish the conditions {\bf BC}, we note that the first is easy because we can take ${\mathcal C}_c^{\infty}(V)$ as a common core for $C$, $D$, and $C+D$. (Note that $D+C$ makes sense as an unbounded self-adjoint operator by the Kato-Rellich theorem.) The second condition holds by virtue of the fact that the Dirac operator and Clifford multiplication by $v$ preserve the Schwartz class of functions. To prove {\bf BC3} we use the Leibniz rule to show that $[D,C]g(v) = g(v)$ for any smooth $g$. Thus, $[D,C]$ extends to a bounded operator.

We may thus invoke the composition formula (Theorem \ref{THM: compform}) to show that the composition $\alpha_*(b) \in E(\C,\C)$ is represented by the operator $D+C$.

The remainder of the proof involves proving that the spectrum of $D+C$ is discrete with one-dimensional kernel and that the asymptotic morphism $\alpha \circ b$ is homotopic to the functional calculus of the projection onto that kernel. This follws from a well-known (but tedious) computation reducing the spectral theory of $D+C$ to that of the one-dimensional harmonic oscillator. Since the functional calculus of a one-dimension projection represents the class of $1$ in $K(\C)$, we are done.
\endproof

\begin{remark}
 The operator $D+C$ appearing in the proof is the so-called Bott-Dirac operator, which has many generalizations including one to infinite-dimensional Euclidean space, see \cite{HKT}. We omit the details here, but the composition formula in Theorem \ref{THM: compform} can be used for (indeed was inspired by) the exact same step of the proof in \cite{HKT}
\end{remark}

\subsection{Perturbation of K-homology Classes}
In this section we will realize many perturbations of elliptic operators as yielding the same K-homology classes as the unpertubed operator. This application is similar to the previous one in terms of the operator theory involved, but seeks to prove a fundamentally different statement about K-theory elements. In this we case, we recognize a composition as such, rather than setting out to compute one.

Let $\xi \to M$ be a smooth, $\Z/2$-graded, Hermitian vector bundle over a smooth compact manifold $M$. Let $D$ be an odd, first-order differential operator acting on smooth sections of $\xi$. The $C^*$-algebra $C(M)$ is represneted on $L^2(\xi)$ as multiplication operators, i.e. there exists a (nondegenerate) $*$-homomorphism $\phi:C(M) \to {\mathcal B}(L^2(\xi))$. The pair $(\phi,D)$ is, in fact, an asymptotic pair (see Section 3 of \cite{Gue1}) and thus defines an element $E'(C(M),{\mathcal K}(L^2(\xi))) \cong K^0(M)$.

We can perturb $D$ by adding a ``zeroth-order potential'' $V$ which is an odd smooth section of the endomorphism bundle ${\rm End}(\xi)$ of $\xi$. Let $D_V = D+V$. Again by the Kato-Rellich theorem, $D_V$ is well-defined as an odd, unbounded, self-adjoint operator on $L^2(\xi)$. It also defines a class in $E'(C(M),{\mathcal K}(L^2(\xi))) \cong K^0(M)$ through the pair $(\phi,D_V)$. We aim to show that these classes are equal.

First we need an easy lemma.

\begin{lemma}\label{LEM: homom}
 Let $A$, $B$ be graded $C^*$-algebras with $B$ stable and let $\phi:A \to B$ a $*$-homomorphism. Let $D$ be any odd, unbounded self-adjoint multiplier of $B$ for which $(\phi,D)$ is an asymptotic pair. Then the classes of $(\phi,D)$ and $(\phi,0)$ are equal in $E'(A,B)$.
\end{lemma}
\proof
By the strict continuity of the functional calculus, for every $a \in A$ and $f \in S$ the family of elements $f(t^{-1}D)\phi(a)$ converges strictly to $f(0)\phi(a)$. This is because $\phi(a)$ actually lies in $B$, not just $M(B)$. Hence, the asymptotic morphisms associated to $(\phi,D)$ and $(\phi,0)$ are asymptotically equivalent.
\endproof

\begin{prop}
 The classes of $(\phi,D)$ and $(\phi,D_V)$ are equal in $E'(C(M),{\mathcal K}(L^2(S)))$.
\end{prop}
\proof
We show that the class of $(\phi,D_V)$ is the composition of the the class of $(\phi,D)$ with the class of the identity in $E'({\mathcal K}(L^2(S)),{\mathcal K}(L^2(S))) \cong K_0({\mathcal  K}(L^2(S))) \cong \Z$.

Note that $V$ defines an odd multiplier of ${\mathcal K}(L^2(S))$. Indeed, $V$ is a bounded multiplier since $M$ is compact. The pair $({\rm id},V)$ is an asymptotic pair determining a class in $E'({\mathcal K}(L^2(S)),{\mathcal K}(L^2(S)))$. From the above lemma we see that the class of $({\rm id},V)$ is the class of the identity homomorphism in $E'({\mathcal K}(L^2(S)),{\mathcal K}(L^2(S)))$.

The push-forward of $D$ by the identity homomorphism is $D$. The operators $D$ and $V$ preserve the space of smooth sections of $S$ and a Leibniz rule calculation shows that $[D,V]$ extends to a bounded operator on $L^2(S)$. Thus, the conditions {\bf BC} are satisfied and we use the composition formula to conclude that the class of $(\phi,D_V)$ is the compostion of the the class $(\phi,D)$ with the class $({\rm id},V)$ and we have the result.
\endproof

\begin{remark}
 A computation like the above proposition works in many instances on noncompact manifolds provided one pays attention to how fast the potential grows at infinity and that the operator $D$ has finite propagation speed. 
\end{remark}

\appendix
\section{Regular Operators and Unbounded Multipliers}\label{SEC: RegOps}
In the following sections we will be using the theory of regular operators on Hilbert modules and, what amounts to the same thing, unbounded multipliers of $C^*$-algebras. We provide the relevant definitions and results here to establish notation. Proofs will be omitted. For a complete treatment of Hilbert modules and the subtle issues regarding unbounded operators on them consult \cite{Lan} and the references therein. For a complete treatment of graded $C^*$-algebras and graded Hilbert modules see \cite{Bla}.

Let $X$ be a right Hilbert $A$-module. We will occasionally use a subscript to emphasize which $C^*$-algebra is acting, e.g. $X_A$. We let ${\cal L}(X)$ denote the $C^*$-algebra of bounded adjointable (with respect to the $A$-valued inner product) operators on $X$ and ${\cal K}(X)$ denote the closure of the algebra of operators of finite $A$-rank, the ``compact'' operators on $X$. 

All $C^*$-algebras and Hilbert modules are assumed separable and $\Z/2\Z$-graded. Recall that a $C^*$-algebra $A$ is called graded if it is equipped with a $*$-automorphism $\gamma \in {\rm Aut}(A)$ such that $\gamma^2 = {\rm id}_A$. Such an automorphism provides  $A$ with a direct sum decomposition $A = A_0 \oplus A_1$, where $A_i$ is the $(-1)^i$-eigenspace of $\gamma$. Elements of $a \in A_i$ are called homogeneous of degree $i$, denoted $\partial a = i$. Alternatively, elements of $A_0$ are called even and elements of $A_1$ are called odd. A grading on $A$ induces a grading on the multiplier algebra $M(A)$ by $T \mapsto \gamma\circ T \circ \gamma$. This grading restricts to the original one on $A$.

A graded Hilbert module is a Hilbert module equipped with a direct sum decomposition $X=X_0 \oplus X_1$. Such a decomposition provides the $C^*$-algebra ${\cal K}(X)$ of compact adjointable operators on $X$ with the grading of conjugation by the unitary ${\rm diag}(1,-1)$, where the matrix is with respect to the decomposition $X = X_0 \oplus X_1$. A grading given by conjugation by a unitary is called an inner grading. 

Important examples of graded $C^*$-algebras for us are:
\begin{itemize}
 \item The algebra $S = C_0(\R)$ of continuous functions on $\R$ which vanish at infinity. The grading is given by even and odd functions, i.e. $\gamma(f)(x) = f(-x)$.
\item The compact operators ${\cal K}(X)$ over a graded Hilbert module $X$. The grading is into diagonal and off-diagonal matrices with respect to the decomposition $X = X_0 \oplus X_1$.
\end{itemize}

A $*$-homomorphism $\phi:A \to B$ between graded $C^*$-algebras is called a {\it graded homomorphism} if $\phi\circ\gamma_A = \gamma_B\circ\phi$. In other words, $\phi(A_i) \subset B_i$. All $*$-homomorphisms will be assumed graded. To each $*$-homomorphism $\phi$ we can associate the opposite homomorphism $\phi^{{\rm opp}}=\phi \circ \gamma$.

\begin{lemma} \label{LEM: gradder} Let $A$ be a graded $C^*$-algebra and let $a$, $b$, and $c$ be homogeneous elements of $A$. Then:
\begin{enumerate}
\item[(i)] $[a,b]  + (-1)^{\partial a \partial b} [b,a] = 0$
\item[(ii)] $[a,bc] = [a,b]c + (-1)^{\partial a \partial b}b[a,c]$
\end{enumerate}
\end{lemma}

\begin{remark} Part (ii) is a graded analogue of the derivation property of ordinary (ungraded) commutators.
\end{remark}

When considering the tensor product of $C^*$-algebras $A$ and $B$, we will always use the maximal graded tensor product, denoted $A \hat{\otimes} B$. See Chapter 14 of Blackadar \cite{Bla} for its construction. The maximal tensor product has the universal mapping property that a pair of $*$-homomorphisms $\phi:A \to C$ and $\psi:B \to C$ that graded commute induces a map $\phi \hat{\otimes} \psi$ on $A \hat{\otimes}B$ given on elementary tensors by $a\hat{\otimes}b \mapsto \phi(a)\psi(b)$.

The {\em graded} commutator of two homogeneous elements $a,b \in A$ is defined as $$[a,b] = ab - (-1)^{\partial a \partial b}ba.$$ All commutators are assumed to be graded commutators. The next lemma collects some easy identities for graded commutators.

We wish to consider closed linear maps on $X$ which are not bounded, but which are still adjointable with respect to the $A$-valued inner product. Let $$D:{\rm dom}(D) \to X$$ be a closed $A$-linear map from a dense right $A$-submodule ${\rm dom}(D)$ of $X$ to $X$. We say that $D$ is symmetric if ${\rm dom}(D) \subseteq {\rm dom}(D^*)$, where $D^*$ is the adjoint of $D$. Recall that $${\rm dom}(D^*) = \{ x \in X | \text{ there exists } y \in X \text{ with } \langle Dz,x \rangle = \langle z,y \rangle \text{ for all } z \in {\rm dom}(D) \}$$
 (c.f. Chapter 9 in \cite{Lan}). We say that $D$ is self-adjoint if it is symmetric and ${\rm dom}(D)={\rm dom}(D^*)$.

There is a further technicality regarding unbounded operators on Hilbert modules that is not present for operators on Hilbert space.

\begin{definition} A self-adjoint operator $D$ on $X$ is called {\bf regular} if the map $1 + D^2$ has dense range.
\end{definition}

Regularity ensures that the graph of $D$ is orthocomplemented with respect to the $A$-valued inner product on $X \times X$, which is implicitly necessary for all of the following results. Of paramount importance for us is that regular self-adjoint operators have a functional calculus (Proposition \ref{PROP: funccalc} below). 

When dealing with unbounded operators on graded Hilbert modules, we restrict ourselves to those operators whose domains are dense, graded submodules. In particular, we assume that if $x \in {\rm dom}(D)$ then its even and odd parts are in ${\rm dom}(D)$. 

\begin{definition}
 A regular, self-adjoint operator $D$ on a graded, right Hilbert $A$-module is called {\bf odd} if $\gamma  D \gamma = -D$. 
\end{definition}
 
Odd multipliers are important because their functional calculus has extra structure, namely, it is a graded $*$-homomorphism.

\begin{prop}[Theorem 10.9 \cite{Lan}] \label{PROP: funccalc} Let $D$ be an odd, regular, self-adjoint operator on $X$. There exists a graded $*$-homomorphism $C_0(\R) \rightarrow {\cal L}(X)$ which sends the resolvent functions $$r_{\pm}(x)=(x \pm i)^{-1}$$ to the operators $$r_{\pm}(D)=(D \pm i)^{-1}.$$ The image of $f \in C_0(\R)$ under this map is denoted $f(D)$. This homomorphism is nondegenerate and thus extends to a strictly continuous $*$-homomorphism $M(C_0(\R)) \cong C_b(\R) \rightarrow {\cal L}({\cal H})$. 
\end{prop}

\begin{remark} Since the functional calculus is a $*$-homomorphism, we have the estimate $$\| f(T) \| \leq \|f\|_{\infty},$$ for every $f \in C_b(\R)$.
\end{remark} 

We now specialize the class of unbounded operators to those which act on $C^*$-algebras. Note that $A$ is a right Hilbert $A$-module over itself.

\begin{definition} An {\bf unbounded self-adjoint multiplier} $D$ of a $C^*$-algebra $A$ is a regular, self-adjoint operator on $A$ viewed as a module over itself. If ${\cal A} \subseteq A$ is a dense right $A$-submodule such that $D$ is the closure of its restriction to ${\cal A}$, then we call ${\cal A}$ a {\bf core} for $D$ and we say that $D$ is {\bf essentially self-adjoint} on ${\cal A}$.
\end{definition}

\begin{remark} Essentially self-adjoint multipliers have a unique self-adjoint extension. If $D$ is essentially self-adjoint on some domain ${\cal A}$, then its resolvents are densely defined, adjointable and bounded. Their closures are the resolvents of the closure of $D$. Thus, essentially self-adjoint multipliers have a well-defined functional calculus. For our purposes, this is enough. All unbounded multipliers will be assumed odd for the reasons mentioned above.
\end{remark}

Recall that a $*$-homomorphism $\phi:A \to M(B)$ is called strict if the image of an approximate unit in $A$ forms a strictly Cauchy net in $M(B)$. Since any multiplier $T \in M(A)$ is the strict limit of $Te_{\lambda}$ for any approximate unit $e_{\lambda}$, we can define an extension of $\phi$ to $M(A)$ by $\phi(T) = \lim_{\lambda \to \infty} Te_{\lambda}$. Furthermore, the image of an approximate unit converges strictly to a projection in $M(B)$ and $\phi$ is a nondegenerate $*$-homomorphism from $A$ to ${\cal L}(X)$, where $X = PB$. 

Just as bounded multipliers of $C^*$-algebras can be pushed-forward by strict $*$-homomorphisms, unbounded multipliers have push-forwards. The purpose here is to collect the relevant results from \cite{Wor} and \cite{WorNap}. The following proposition was stated in \cite{Wor} for push-forwards by nondegenerate $*$-homomorphisms. We state it here for strict homomorphisms. There is no added difficulty to the proof.

\begin{prop}[Theorem 1.2 \cite{Wor}] \label{PROP: pushforward} Let $\phi:A \rightarrow M(B)$ be a strict $*$-homomorphism and let $D$ be an unbounded multiplier of $A$. Then there exists an unbounded multiplier of $B$ denoted $\phi(D)$ which is defined on the core ${\rm dom}(\phi(D)) = \phi({\rm dom}(D))\cdot B$ and is given by the formula $$\phi(D) (\phi(a)b) = \phi(Da)b.$$ Furthermore, if $\psi:B \rightarrow M(C)$ is another strict $*$-homomorphism, then $\psi(\phi(T)) = (\psi \circ \phi)(T)$.
\end{prop}

 \begin{remark} The formula for $\phi(D)$ is exactly the same as for the push-forward of a bounded multiplier, but in this case it is more difficult to show that it is well-defined. 
 \end{remark}

The set of unbounded multipliers of $C_0(\R)$ is $C(\R)$, where a continuous (not necessarily bounded) function acts by multiplication (c.f. page 117 of \cite{Lan}). Such operators are all defined on the core $C_c(\R)$. Since the functional calculus is a nondegenerate (hence strict) $*$-homomorphism, we have the following corollary of Proposition \ref{PROP: pushforward} which gives a functional calculus for unbounded continuous functions.

\begin{corollary}[Theorem 10.9 \cite{Lan}]
\label{COR: unboundedfunccalc}
 Let $D$ be an unbounded self-adjoint multiplier of $A$. There is an adjoint preserving map from $C(\R)$ to the set of unbounded multipliers of $A$. The image of $f \in C(\R)$ is denoted $f(D)$ and is essentially self-adjoint on the domain $C_c(\R)\cdot A$.
\end{corollary}
 
 There is another important corollary to  Proposition \ref{PROP: pushforward} which can be stated informally as, ``functional calculus commutes with push-forwards".
 
 \begin{corollary}[Theorem 1.2 \cite{Wor}]
 Let $A$ and $B$ be $C^*$-algebras, $\phi : A \rightarrow M(B)$ a strict $*$-homomorphism, and $D$ an unbounded, self-adjoint multiplier of $A$ with domain $\mathcal{A}$. For every $f \in C(\R)$, $f(\phi(D)) = \phi(f(D))$.
\end{corollary}
 
 As a final note on push-forwards, we can relax the requirement in Proposition \ref{PROP: pushforward} that $D$ be self-adjoint to requiring that $D$ be essentially self-adjoint. In this case the push-forward of an essentially self-adjoint multiplier is essentially self-adjoint.

To end this subsection we would like to state a useful relationship between regular operators on a Hilbert module and unbounded multipliers of the the compacts over that Hilbert module. This observation is well-known, but to the author's knowledge it was first put in print in \cite{Pal}.

\begin{prop}[Section 3 \cite{Pal}]  \label{PROP: hilbmodop} Let $X$ be a right Hilbert $A$-module. There is a one-to-one correspondence between regular self-adjoint operators on $X$ and unbounded self-adjoint multipliers of ${\cal K}(X)$. If $D$ is a regular self-adjoint operator on $X$, then the associated unbounded self-adjoint multiplier has domain $\{T \in {\cal K}(X): {\rm ran}(T) \subseteq {\rm dom}(D)$ and $DT \in {\cal  K}(X)\}$ and is given by left multiplication on this domain.
\end{prop}

\begin{remark}
 The reverse direction of the correspondence is given by identifying $X$ with the set of rank one projections in ${\cal K}(X)$, but our main interest is in the forward direction. 
\end{remark}

\section{An Exponential Computation}
In this appendix we will provide a detailed computation that is required in the proof of the composition formula. Stated briefly, we need an identity  relating the exponential of a sum to the separate exponentials of the summands. The sticking points are that we are working with unbounded multipliers and that we are considering graded commutators. Furthermore, we only assume that the multipliers asymptotically commute. The morass of detail here warrants careful consideration.

\begin{lemma} \label{LEM: exp1} Let $C$ be a $C^*$-algebra. Suppose that  $x,y \in C$ such that $\|y\| \leq \|x\|$. Then $$\|e^{x+y}-e^x\| \leq \|y\| e^{2\|x\|}.$$
\end{lemma}
\proof
Using the (absolutely convergent) power series representations of $e^{x+y}$ and $e^x$, we see that $$e^{x+y} - e^x = \sum_{n=0}^{\infty} \frac{1}{n!}[ (x+y)^n - x^n ].$$ Since $x$ and $y$ do not necessarily commute, we must take order into account when expanding the binomial $(x+y)^n$. The first term, $x^n$, gets canceled, but the rest remain and each has at least one factor of $y$ in it. Using the triangle inequality
\begin{align*}
\Big \|\sum_{n=0}^{\infty} \frac{1}{n!}[ (x+y)^n - x^n ] \Big \| &\leq \sum_{n=0}^{\infty} \frac{1}{n!}\|(x+y)^n - x^n\| \\ &= \sum_{n=0}^{\infty} \frac{1}{n!} \|y\| \sum_{j=0}^{n-1} {n \choose j} \|x\|^{j} \|y\|^{n-j-1} \\ &\leq \sum_{n=0}^{\infty} \frac{1}{n!} \|y\| \sum_{j=0}^{n-1} {n \choose j} \|x\|^{j} \|x\|^{n-j-1} \\ &\leq \|y\| e^{2\|x\|}.
\end{align*}
This complete the proof.
\endproof

\begin{remark} In the above proof, we only invoke the binomial theorem after we have taken the norms of the summands. This is what allows us to pull $\|y\|$ out of the sum. 
\end{remark}

\begin{lemma} \label{LEM: exp2} Let $x_t$ and $y_t$, $t \in [1,\infty)$ be bounded continuous paths of even elements in a graded $C^*$-algebra $C$ such that $[x_t,y_t] \rightarrow 0$ as $t \rightarrow \infty$. Then $$\lim_{t\rightarrow \infty} \|e^{x_t+y_t}- e^{x_t}e^{y_t}\| = 0.$$
\end{lemma}
\proof
The product $e^{x_t}e^{y_t}$ has power series $$e^{x_t}e^{y_t} = \sum_{n=0}^{\infty}\sum_{j+k=n} \frac{1}{j!k!} x_t^jy_t^k.$$ We want to take the difference $e^{x_t+y_t}- e^{x_t}e^{y_t}$, but the series $$e^{x_t+y_t} = \sum_{n=0}^{\infty} \frac{1}{n!} (x_t+y_t)^n, $$ has binomials in the noncommuting elements $x_t$ and $y_t$. When we expand $(x_t+y_t)^n$ we get monomials of the form $x_t^{i_1}y_t^{i_2} \dots y_t^{i_l}$, where $i_1,\dots,i_k$ are exponents which add up to $n$. Our strategy is to rearrange these monomials so that all the $x_t$'s are on the left and $y_t$'s are on the right. In order to do this, we need to add terms containing commutators, e.g. $x_t y_t x_t = x_t^2 y_t - x_t [x_t,y_t]$. Notice that for every swap we add a term with one commutator. The maximum number of swaps needed for a degree $n$ term is $n^2/4$. Indeed, for a term with $j$ $x_t$'s and $(n-j)$ $y_t$'s, one must make at most $j(n-j)$ transpositions. This expression is maximized when $j=n/2$ and for a term of the form $x_t^{n/2}y_t^{n/2}$ one does need $n^2/4$ swaps. 

In the expansion of $(x_t+y_t)^n$ there are $n \choose j$ terms with $j$ $x_t$'s and $(n-j)$ $y_t$'s. After rearrangement, there are ${n \choose j}$ terms of the form $x_t^j y_t^{n-j}$. Let  $M > \sup_{t\in [1,\infty)} \{\|x_t\|,\|y_t\| \}$. Computing, we obtain: 
\begin{align*}
\Big \|\sum_{n=0}^{\infty}[\frac{1}{n!}(x_t+y_t)^n - \sum_{j+k=n} \frac{1}{j!k!} x_t^jy_t^k] \Big \| &\leq \sum_{n=0}^{\infty}\| \frac{1}{n!}(x_t+y_t)^n - \sum_{j+k=n} \frac{1}{j!k!} x_t^jy_t^k \| \\ &\leq \sum_{n=0}^{\infty} \sum_{j=0}^n \frac{1}{n!} {n \choose j} \frac{n^2}{4} \|[x_t,y_t]\| M^{n-2}.
\end{align*}
The inequality on the last line comes from the observation that after reordering the monomials in the binomial expansion, the terms of the form $x_t^j y_t^{n-j}$ all cancel (there are $n \choose j$ of them in both sums). The remaining terms all have one commutator in them and $(n-2)$ $x_t$ or $y_t$ factors. Estimating the sum 
\begin{align*}
\sum_{j=0}^n \frac{1}{n!} {n \choose j} \frac{n^2}{4} \|[x_t,y_t]\| M^{n-2} &\leq (n+1)\sup_{j} \{ \frac{1}{j! (n-j)!} \} \frac{n^2}{4} \|[x_t,y_t]\| M^{n-2} \\ &\leq (n+1) \frac{1}{(n/2)!^2}  \frac{n^2}{4} \|[x_t,y_t]\| M^{n-2}.
\end{align*}
The ratio test shows that this series converges. The factor $\|[x_t,y_t]\|$ goes to zero as $t \rightarrow \infty$.
\endproof 

Now we can provide a proof of Corollary \ref{COR: fclemma}.  Recall that we let $D_t = t^{-1}D$ and $D_{t,N} = (D_t)_N$.

\begin{corollary} Let $D$ and $D'$ be odd, essentially self-adjoint multipliers of a $C^*$-algebra $C$ which have bounded commutator. Then $$\lim_{t \rightarrow \infty} \left( e^{-(D_t+D_t')^2} - e^{-D_t^2}e^{-D_t'^2} \right) = 0.$$
\end{corollary}
\proof 
If both of the multipliers $D$ and $D'$ are bounded, then this proposition is easy. We will assume, then, that at least one of $D$ or $D'$ is unbounded.

The proof will require several applications of the ``add zero" trick. We want to show that the difference $$e^{-(D_t+D_t')^2} - e^{-D_t^2}e^{-D_t'^2}$$ goes to zero as $t \rightarrow \infty$. Adding and subtracting $e^{-(D_{t,N}+D_{t,N}')^2}$, we know by Lemma \ref{PROP: techlemma} that $$\lim_{N\rightarrow \infty} \limsup_{t \rightarrow \infty} \|e^{-(D_t+D_t')^2} - e^{-(D_{t,N}+D_{t,N}')^2}\| = 0.$$ Thus, we aim to show that $$\lim_{N\rightarrow \infty} \limsup_{t \rightarrow \infty} \|e^{-(D_{t,N}+D_{t,N}')^2} -  e^{-D_t^2}e^{-D_t'^2}\| = 0.$$ 

Since $D_N$ and $D_N'$ are both odd, bounded multpliers, $$(D_N+D_N')^2 = D_N^2+D_N'^2+[D_N,D_N'].$$ From Lemma \ref{LEM: commbound} we know that $[D_N,D_N']$ is uniformly bounded in $N$. On the other hand, the norm of $D_N^2+D_N'^2$ is unbounded as $N\rightarrow \infty$, because as positive operators $$D_N^2+D_N'^2 \geq D_N^2$$ and at least one of $D$ or $D'$ is unbounded. Thus, there is an $N>0$ for which $\|D_N^2+D_N'^2\| \geq \|[D_N,D_N']\|$ and we can use the estimate of Lemma \ref{LEM: exp1} to get $$\| e^{-(D_{t,N}+D_{t,N}')^2} - e^{-(D_{t,N}^2+D_{t,N}'^2)} \| \leq \|[D_{t,N},D_{t,N}']\|e^{2\|D_N^2+D_N'^2\|}.$$ For fixed $N$, this difference goes to zero as $t \rightarrow \infty$. Thus, after adding and subtracting $e^{-D_{t,N}^2+D_{t,N}'^2}$, it suffices to show that $$\lim_{N\rightarrow \infty} \limsup_{t \rightarrow \infty} \|e^{-(D_{t,N}^2+D_{t,N}'^2)} -  e^{-D_t^2}e^{-D_t'^2}\| = 0.$$

Letting $x_t=-D_{t,N}^2$ and $y_t=-D_{t,N}'^2$, we see that $x_t$ and $y_t$ are bounded, continuous paths of even elements in $M(A)$. Using Lemma \ref{LEM: gradder}, we see that $\|[x_t,y_t]\|\rightarrow 0$ as $t\rightarrow \infty$. From Lemma \ref{LEM: exp2} we know that $$\|e^{x_t+y_t}- e^{x_t}e^{y_t}\|  \rightarrow 0$$ as $t \rightarrow \infty$. Hence after adding and subtracting $e^{-D_{t,N}^2}e^{-D_{t,N}'^2}$, it suffices to show that $$\lim_{N\rightarrow \infty} \limsup_{t \rightarrow \infty} \|e^{-D_{t,N}^2}e^{-D_{t,N}'^2} - e^{-D_t^2}e^{-D_t'^2}\| = 0,$$ but this is an application of Lemma \ref{PROP: techlemma}. This completes the proof.
\endproof 

\begin{remark}
 The above corollary only provides the first part of \ref{COR: fclemma}. The second part is simple an application of \ref{COR: asycomm2} after terms of the form $t^{-1}De^{t^{-2}D^2}$ are approximated by $t^{-1}D_Ne^{t^{-2}D^2}$.
\end{remark}

If $[D,D']=0$, then an inspection of the above proofs shows that, because $[D_N,D_N'] = 0$, all of the above combinatorial arguments go through without having to bound error terms. Corollary \ref{COR: fclemma2} follows immediately.

\bibliographystyle{amsplain}
\bibliography{mattref}                                                                                                              
                                                                                                                                                                                                
% Set the ending of a LaTeX document                                                                                                                                                            
\end{document}